\documentclass{article}
\usepackage[left=2.5cm,right=2.5cm,top=2.5cm,bottom=2.5cm]{geometry}
\usepackage{amssymb,amsmath,amsthm}
\usepackage{paralist}
\usepackage{graphics} 
\usepackage{epsfig} 
\usepackage{graphicx}
\usepackage{epstopdf}
\usepackage[colorlinks=true]{hyperref}

\newtheorem{remark}{Remark}

\newcommand{\R}{\mathbb{R}}

\newcommand{\fref}[1]{Fig.~\ref{#1}}
\newcommand{\Fref}[1]{Figure~\ref{#1}}
\title{Computing connecting orbits to infinity \\
associated with a homoclinic flip bifurcation}

\author{Andrus Giraldo\footnotemark[1], Bernd Krauskopf\footnotemark[1]  and Hinke M. Osinga\footnotemark[1]}

\date{}

\begin{document}

\renewcommand{\thefootnote}{\fnsymbol{footnote}}
\footnotetext[1]{Department of Mathematics, The University of
  Auckland, Private Bag 92019, Auckland 1142, New Zealand
  (\href{mailto:a.giraldo@auckland.ac.nz}{a.giraldo@auckland.ac.nz},
  \href{mailto:b.krauskopf@auckland.ac.nz}{b.krauskopf@auckland.ac.nz}, \href{mailto:h.osinga@auckland.ac.nz}{h.osinga@auckland.ac.nz})}
\renewcommand{\thefootnote}{\arabic{footnote}}

\maketitle

\begin{abstract}
We consider the bifurcation diagram in a suitable parameter plane of a quadratic vector field in $\R^3$ that features a homoclinic flip bifurcation of the most complicated type. This codimension-two bifurcation is characterized by a change of orientability of associated two-dimensional manifolds and generates infinite families of secondary bifurcations. We show that curves of secondary $n$-homoclinic bifurcations accumulate on a curve of a heteroclinic bifurcation involving infinity. 

We present an adaptation of the technique known as Lin's method that enables us to compute such connecting orbits to infinity. We first perform a weighted directional compactification of $\R^3$ with a subsequent blow-up of a non-hyperbolic saddle at infinity. We then set up boundary-value problems for two orbit segments from and to a common two-dimensional section: the first is to a finite saddle in the regular coordinates, and the second is from the vicinity of the saddle at infinity in the blown-up chart. The so-called Lin gap along a fixed one-dimensional direction in the section is then brought to zero by continuation. Once a connecting orbit has been found in this way, its locus can be traced out as a curve in a parameter plane.
\end{abstract}

\section{Introduction}

Homoclinic flip bifurcations are bifurcations of codimension two that occur in families of continuous-time dynamical systems, given by ODEs or vector fields, whose phase space is of dimension at least three. This type of bifurcation of a homoclinic orbit to a real hyperbolic saddle---a special trajectory that converges both in forward and backward time to the saddle equilibrium---occurs when a stable or unstable manifold transitions, when followed along the homoclinic orbit, from being orientable to being non-orientable, or vice versa. While such a change of orientability may occur in higher-dimensional phase spaces, the characterization of homoclinic flip bifurcations and their unfoldings have been studied in detail mostly for the lowest-dimensional case of a three-dimensional systems, both from a theoretical~\cite{Deng1993, HKK1994, HKN2001, HK2000, KKO1993, OKC2001} and a numerical point of view~\cite{AKO2013, CK1994, And1, And2018, san3}. 

In three-dimensions, which is the case we also consider here, the orientability of a homoclinic orbit is determined by the orientablility of the two-dimensional (un)stable manifold. The saddle equilibrium is assumed to be a hyperbolic saddle, meaning that it has one or two stable and two or one unstable eigenvalues, respectively. In the case of one stable eigenvalue, which we encounter in the example vector field below, its stable manifold is one dimensional, that is, a curve consisting of two trajectories that converge to the saddle in forward time; its unstable manifold is two dimensional, that is, a surface formed by all trajectories that converge to the saddle in backward time. Generically, this surface, when followed locally along the homoclinic orbit in backward time to the equilibrium, closes up along the one-dimensional strong unstable manifold, which is tangent to the strongest unstable eigendirection of the saddle, to form either a cylinder in the orientable case, or a M{\"o}bius strip in the non-orientable case. The orientability of the homoclinic orbit, that is, of the two-dimensional unstable manifold (in this case), can change in three different ways: \\[-4mm]
\begin{list}{\arabic{enumi}.}{\usecounter{enumi}
    \setlength{\leftmargin}{5mm} \setlength{\labelwidth}{5mm} \setlength{\itemsep}{1mm}}
\item the two unstable eigenvalues become complex conjugate and the equilibrium turns into a saddle-focus;
\item orbit flip: the one-dimensional stable manifold returns (in backward time) to the equilibrium tangent to the strong unstable eigendirection (instead of the weakest unstable eigendirection);
\item inclination flip: the two-dimensional unstable manifold when followed along the homoclinic orbit is tangent to the plane spanned by the stable and weak unstable eigendirections (instead of the plane spanned by the stable and strong unstable eigendirections). \\[-3mm]
\end{list}
The first case, when the saddle equilibrium has a double leading eigenvalue, is known as a Belyakov bifurcation~\cite{bel80, bel84} and a numerical study of its simplest unfolding was performed in~\cite{CK1994, KFR2001}; see also~\cite{AKO2014}. Both the orbit flip and inclination flip bifurcations have similar unfoldings, which are split into three different generic cases, referred to as $\mathbf{A}$, $\mathbf{B}$ and $\mathbf{C}$, depending on the eigenvalues of the saddle equilibrium; see~\cite{HKK1994, HKN2001, HK2000, san3} for the actual eigenvalue conditions. In a two-parameter unfolding, case~$\mathbf{C}$, which is the most complicated one, gives rise to infinitely many curves of secondary bifurcations, including saddle-node, period-doubling, and $n$-homoclinic bifurcations (which involve so-called $n$-homoclinic orbits that make $n-1$ close passes of the equilibrium before returning to it). Moreover, there are two different unfoldings with quite different arrangements of the associated secondary bifurcations, called inward twist $\mathbf{C}_{\rm in}$ and outward twist $\mathbf{C}_{\rm out}$; which of the two unfoldings occurs is determined by global geometric properties of the two-dimensional manifold~\cite{Deng1993, HKK1994}. 

We previously conducted numerical studies of the unfoldings of the different homoclinic flip bifurcations, with a particular focus on clarifying changes of two-dimensional global invariant manifolds~\cite{AKO2013, And1, And2018}. To this end, we studied a model vector field developed by Sandstede~\cite{san4, Sandstede1997}: a system of three ordinary differential equations with eight parameters, which contains all three cases $\mathbf{A}$, $\mathbf{B}$ and $\mathbf{C}$ of both orbit flip and inclination flip bifurcations for suitable choices of the parameters; the underlying homoclinic orbit is always to the saddle located at the origin. However, all unfoldings of case~$\mathbf{C}$ in Sandstede's model are outward twisted for the inclination flip~\cite{Sandstede1997}. Furthermore, we considered the case of the orbit flip in Sandstede's model and did not find a parameter regime where the unfolding of case~$\mathbf{C}$ is inward twisted. In fact, no explicit example of a vector field with the inward-twisted case $\mathbf{C}_{\rm in}$ of a flip bifurcation was known.

This has changed very recently, when Algaba, Dom\'{\i}nguez-Moreno, Merino, and Rodr\'{\i}guez-Luis~\cite{Alg2019} found an example of a three-dimensional quadratic system with an inward-twisted homoclinic flip bifurcation. More precisely, they presented the system
\begin{equation}
\label{eq:AlgabaEq}
  \left\{ \begin{array}{rcrcl}
    \dot x_1 &=&    a \ &+& x_2 \, x_3, \\[1mm]
    \dot x_2 &=& -x_2 &+&  x_1^2, \\[1mm]
    \dot x_3 &=&   b \ &-& 4  x_1
  \end{array} \right.
\end{equation} 
and showed that it exhibits a codimension-two homoclinic orbit flip bifurcation of type $\mathbf{C}_{\rm in}$ of a saddle $p$ when $a \approx -1.20338$ and $b \approx 1.89616$. This was achieved by identifying the orbit flip homoclinic bifurcation numerically in a parameter regime where the eigenvalue condition at $p$ of case $\mathbf{C}$ is satisfied, and then computing a sufficient number of secondary bifurcation curves emanating from this codimension-two point to show that it unfolds as case $\mathbf{C}_{\rm in}$. Note that the homoclinic orbit is not to the origin but to the equilibrium  $\mathbf{p} = (x_1, x_2, x_3) = (b / 4,\, b^2 / 16,\, -16 \, a / b^2)$, which exists provided $b \neq 0$.

Algaba \emph{et al.}~\cite{Alg2019} studied the local bifurcation structure near $\mathbf{C}_{\rm in}$ in quite some detail. We are interested here in how the unfolding of $\mathbf{C}_{\rm in}$ is embedded more globally in an overall bifurcation diagram. An interesting aspect of system~\eqref{eq:AlgabaEq} is that it has only one finite equilibrium, the origin that is involved in the homoclinic bifurcation. By contrast, in Sandstede's model there exists a second equilibrium, and we found that it is responsible for additional global bifurcations in the overall bifurcation diagram, including connecting orbits to the origin~\cite{AKO2013, And1, And2018}.

In this paper, we focus on a particular global feature of system~\eqref{eq:AlgabaEq}, namely connecting orbits to a second equilibrium $\mathbf{q}_\infty$ that, intriguingly, is located at infinity. More specifically, we study the bifurcation diagram near $\mathbf{C}_{\rm in}$ in a suitable two-parameter plane and show that it features curves of $n$-homoclinic bifurcation that emanate from $\mathbf{C}_{\rm in}$. We find that these curves accumulate, as $n$ increases, on a curve of heteroclinic bifurcations involving infinity, given by the existence of a connecting orbit of codimension-one from the finite equilibrium $\mathbf{p}$ to the equilibrium $\mathbf{q}_\infty$ at infinity. Hence, the bifurcation diagram near the codimension-two orbit flip point $\mathbf{C}_{\rm in}$ in the quadratic system~\eqref{eq:AlgabaEq} features global connecting orbits to infinity. 

To address the challenge of finding such connecting orbits to infinity, we adapt the numerical technique from~\cite{KraRie1}, referred to as Lin's method, for computing connecting orbits between finite objects. More precisely, we modify system~\eqref{eq:AlgabaEq} by translating the equilibrium $\mathbf{p}$ to the origin $\mathbf{0}$ and by introducing a third parameter that helps separate the very closely spaced bifurcations. For the transformed system, we perform a weighted directional compactification of phase space to study the behavior at infinity. The analysis at infinity involves an additional blow-up transformation to understand the behavior of solutions approaching $\mathbf{q}_\infty$ in backward time; these are bounded in the blow-up chart by a two-dimensional surface related to a specific periodic orbit at infinity. To set up Lin's method, we choose a section $\Sigma$ that is well defined in the original coordinates as well as the blow-up chart near infinity. We then consider and compute two orbit segments, from this periodic orbit surrounding $\mathbf{q}_\infty$ to $\Sigma$ and from $\Sigma$ to $\mathbf{0}$, such that their end points in $\Sigma$ lie in the so-called Lin space. In this way, we obtain a well-defined and computable test function, which is zero exactly at the parameter values where there exists a connecting orbit to $\mathbf{q}_\infty$. All our computations are performed via the continuation of solutions to suitable two-point boundary value problems with the pseudo-arclength continuation package \textsc{Auto}~\cite{Doe1, Doe2} and the homoclinic continuation toolbox \textsc{HomCont}~\cite{san2}. 

This paper is organized as follows. In the next section, we introduce the transformed system with a homoclinic orbit to the origin. Furthermore, we identify the codimension-two point $\mathbf{C}_{\rm in}$ and present a bifurcation diagram in two parameters that suggests the need for the analysis of the dynamics at infinity. Section~\ref{sec:compact} presents the compactification and the blow-up analysis in different charts at infinity. We use these results in Section~\ref{sec:LinAtInfty} to set up Lin's methods by defining a suitable boundary value problem to compute the boundary of the existence of connecting orbit from $\mathbf{0}$ to the equilibrium $\mathbf{q}_\infty$ at infinity. Section~\ref{sec:LinPeriodicInfty} then explains how this set-up can also be used to find connecting orbits from a saddle periodic orbit to infinity. In the final  Section~\ref{sec:conclusions} we draw conclusions and point to some directions for further research.

\section{Codimension-two orbit flip bifurcation of inward-twisted type $\mathbf{C}_{\rm in}$}
A homoclinic flip bifurcation of case $\mathbf{C}$ is the global bifurcation of the lowest codimension that involves a real saddle equilibrium (its eigenvalues relevant to the bifurcation are all real) and gives rise to chaotic dynamics. While its complete unfolding is not fully understood, a lot is known about the dynamics nearby~\cite{HKK1994, HKN2001, HK2000, KKO1993, OKC2001, And2018, san3}. It has been proven that there exists a nearby parameter region with Smale-horseshoe dynamics, and this means that infinitely many saddle periodic orbits are created near this codimension-two point. The precise way in which this occurs is organized by cascades of period-doubling and saddle-node bifucations as well as cascades of $n$-homoclinic bifurcations; these infinitely many different bifurcations occur arbitrarily close in parameter space to the homoclinic flip bifurcation point. The difference between the two cases $\mathbf{C}_{\rm our}$ and $\mathbf{C}_{\rm in}$ lies in the positions of these cascades relative to the primary homoclinic orbit that undergoes the flip bifurcation.

Algaba \emph{et al.}~\cite{Alg2019} identified an orbit flip bifurcation of system~\eqref{eq:AlgabaEq}, and computed and presented several bifurcation curves in the $(a, b)$-parameter plane to show that the bifurcation diagram is that of inward-twisted type $\mathbf{C}_{\rm in}$. Unfortunately, the bifurcations for system~\eqref{eq:AlgabaEq} occur extremely close together and it is not easy to distinguish them. Furthermore, the multi-loop periodic orbits that are created in the $n$-homoclinic bifurcations come very close to the saddle equilibrium and do not extend far in phase space. In a bid to ameliorate this, we move the unique equilibrium $\mathbf{p}$ of~\eqref{eq:AlgabaEq} to the origin and introduce a third parameter to obtain the system
\begin{equation}
\label{eq:WorkEquation}
  \left\{ \begin{array}{rcrcccl}
            \dot{x} &=& \alpha \, y &+& \gamma \, z &+& y \, z, \\
            \dot{y} &=& \beta \, x &-& y &+& x^2, \\
            \dot{z} &=& -4 \, x.&&&&
          \end{array} \right.
\end{equation} 
Here, the new variables are given by $(x, y, z) = (x_1- b / 4,\, x_2- b^2 / 16,\, x_3 + 16 \, a / b^2)$, and we consider system~\eqref{eq:WorkEquation} as a new system with three independent parameters. Note that the previous system~\eqref{eq:AlgabaEq} is recovered for the special choice of the new  parameters 
$\alpha = -16 a / b^2, \, \beta = b / 2 \text{ and } \gamma = b^2 / 16.$
The advantage of having the specific parameter $\gamma$ is that it allows us to improve the separation of bifurcating periodic orbits from $\mathbf{0}$, the only equilibrium of system~\eqref{eq:WorkEquation}. For the purpose of this paper, the parameters $\alpha$ and $\beta$ are allowed to vary as the unfolding parameters of the orbit flip bifurcation, while we fix $\gamma = 0.5$ throughout our investigation.

System~\eqref{eq:WorkEquation} is our object of study. To find the orbit flip bifurcation in the $(\alpha,\beta)$-plane for $\gamma = 0.5$, we start from the parameter values corresponding to those reported in~\cite{Alg2019} and continue the (primary) homoclinic bifurcation to $\gamma = 0.5$. Next, we continue the locus of the homoclinic bifurcation as a curve in the $(\alpha, \beta)$-plane while keeping $\gamma = 0.5$ fixed throughout all subsequent computations. On the curve of homoclinic bifurcations, we detect the orbit flip point, which we denote $\mathbf{C}_{\rm in}$, at $(\alpha, \beta) \approx (5.3573, 2.19173)$. At this parameter point the origin $\mathbf{0}$ has eigenvalues $\lambda^s \approx -3.7444$, $\lambda^u \approx 0.2108$, and $\lambda^{uu} \approx 2.5335$. Hence, at $\mathbf{C}_{\rm in}$, and also nearby, the point $\mathbf{0}$ is a hyperbolic saddle with a one-dimensional stable manifold $W^s(\mathbf{0})$ and a two-dimensional unstable manifold $W^u(\mathbf{0})$. Moreover, the condition $\mid\! -\lambda^{uu} \!\mid < -\lambda^s$ on the eigenvalues for an orbit flip of case~$\mathbf{C}$ is indeed satisfied at $\mathbf{C}_{\rm in}$~\cite{san3, san4}.  

%
\begin{figure}[t!]
\begin{center}
  \includegraphics{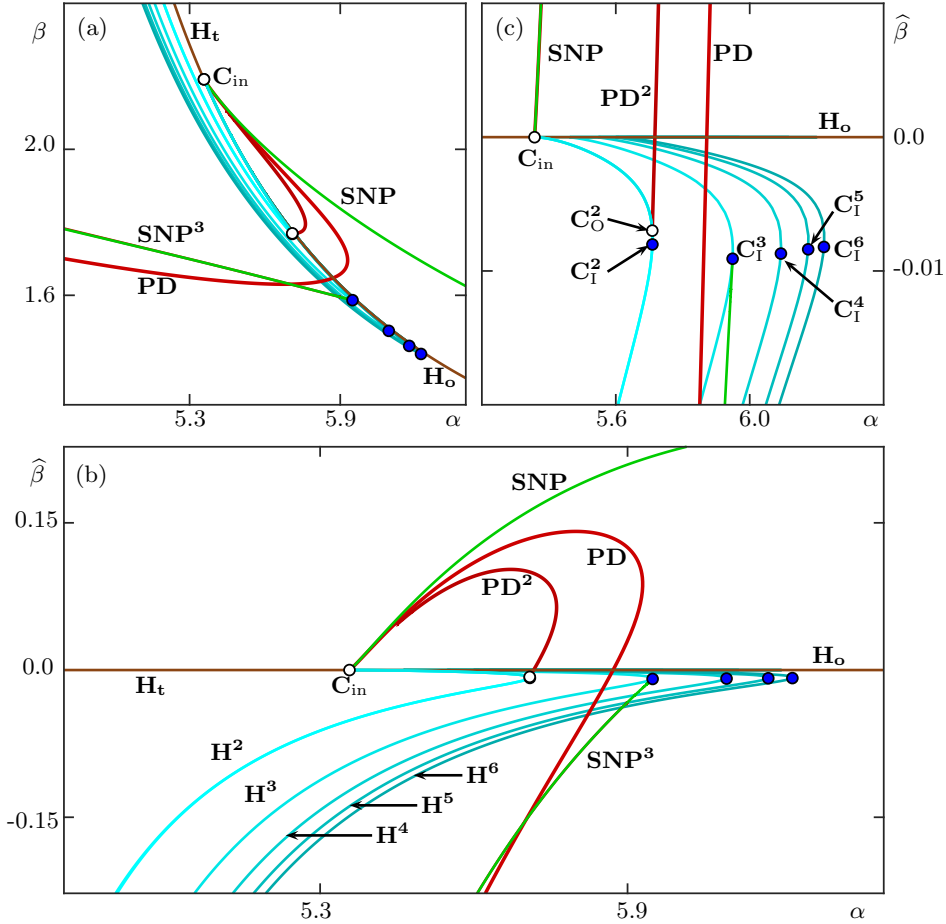}
  \caption{\label{fig:bifDiagram}
    Bifurcation diagram of system~\eqref{eq:WorkEquation} showing: the curve of primary homoclinic bifurcation (brown), along which the homoclinic orbit changes at $\mathbf{C}_{\rm in}$ from being orientable along $\mathbf{H_o}$ to being non-orientable along $\mathbf{H_t}$; curves $\mathbf{SNP}$ and $\mathbf{SNP^3}$ (green) of saddle-node bifurcation of periodic orbits; the first two curves $\mathbf{PD}$ and $\mathbf{PD^2}$ (red) of a cascade of period-doubling bifurcations; and the curves $\mathbf{H^n}$ (increasingly darker shades of cyan) of $n$-homoclinic bifurcations for $n = 2, 3, 4, 5$, and $6$. On $\mathbf{H^n}$ there are points $\mathbf{C^n_{\rm O}}$ of orbit flip bifurcations (blue dots) and on $\mathbf{H^2}$ there is a point $\mathbf{C^2_{\rm I}}$ of inclination flip bifurcation (open dot). Panel~(a) shows the $(\alpha, \beta)$-plane, while panel~(b) shows the $(\alpha, \hat{\beta})$-plane, where $\hat{\beta}$ is the distance in the $\beta$-coordinate from the curve $\mathbf{H_{o/t}}$ of primary homoclinic bifurcation, which is now at $\hat{\beta} = 0$ (brown horizontal line). Panel~(c) is an enlargement of the $(\alpha, \hat{\beta})$-plane near $\mathbf{C}_{\rm in}$.}
\end{center}
\end{figure}
%

\Fref{fig:bifDiagram} shows the partial bifurcation diagram for system~\eqref{eq:WorkEquation}, which provides the numerical evidence that we are indeed dealing with an orbit flip of inward-twisted type $\mathbf{C}_{\rm in}$. The curve of (primary) homoclinic bifurcation in the $(\alpha, \beta)$-plane is separated by the orbit flip point $\mathbf{C}_{\rm in}$ into a branch $\mathbf{H_o}$ of orientable and a branch $\mathbf{H_t}$ of non-orientable or twisted homoclinic bifurcation. Subsequently, we found and continued other bifurcation curves emanating from $\mathbf{C}_{\rm in}$, namely, curves $\mathbf{SNP}$ of saddle-node bifurcation of periodic orbits (green), $\mathbf{PD}$ and $\mathbf{PD^2}$ of period-doubling bifurcation, and $\mathbf{H^n}$ of $n$-homoclinic bifurcation for $n = 2,\dots, 6$. \Fref{fig:bifDiagram}(a) shows the bifurcation diagram in the $(\alpha,\beta)$-plane of~\eqref{eq:WorkEquation}. Because the different bifurcation curves are still a bit hard to distinguish in the $(\alpha,\beta)$-plane, panel~(b) shows them relative to the curve $\mathbf{H_{o/t}}$ of primary homoclinic bifurcation. More specifically, we show the $(\alpha, \hat{\beta})$-plane, where $\hat{\beta}$ represents the distance to $\mathbf{H_{o/t}}$ with respect to the $\beta$-coordinate. Hence, the curve $\mathbf{H_{o/t}}$ is now the $\alpha$-axis where $\hat{\beta} = 0$. \Fref{fig:bifDiagram}(b) illustrates that all bifurcation curves emanate from the point $\mathbf{C}_{\rm in}$ on the side of $\mathbf{H_{o}}$; in particular, the curves $\mathbf{H^n}$ of $n$-homoclinic bifurcation are tangent to $\mathbf{H_{o}}$ near $\mathbf{C}_{\rm in}$, as can be seen in panel~(c). Moreover, the curve $\mathbf{SNP}$, as well as the first two curves $\mathbf{PD}$ and $\mathbf{PD^2}$ of a cascade of period-doubling bifurcations lie on one side of $\mathbf{H_{o/t}}$, while the curves $\mathbf{H^n}$ lie on the other side. These are all characteristic features that distinguish the inward twist from the outward twist~\cite{And2018, san4, Sandstede1997}. Hence, we conclude that the codimension-two point $\mathbf{C}_{\rm in}$ of~\eqref{eq:WorkEquation} for $\gamma=0.5$ is of the same inward-twisted type as that of \eqref{eq:AlgabaEq} found in~\cite{Alg2019}.

%
\begin{figure}[t!]
\begin{center}
  \includegraphics{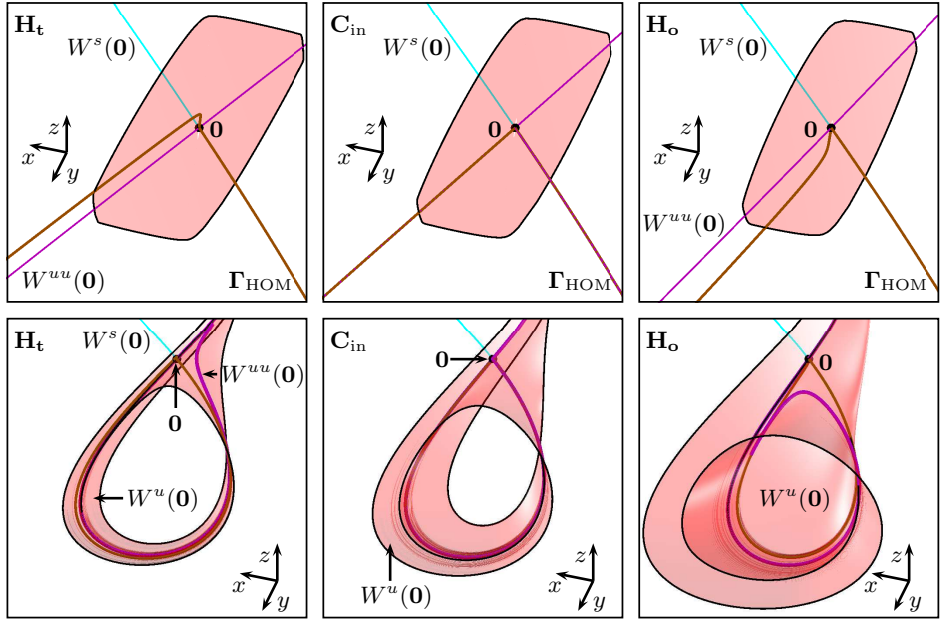}
  \caption{\label{fig:flipTrans}
    Phase portraits of system~\eqref{eq:WorkEquation} along $\mathbf{H_{t}}$, at $\mathbf{C}_{\rm in}$ and along $\mathbf{H_{o}}$ with enlargements near the saddle $\mathbf{0}$ (top row). Shown are the saddle $\mathbf{0}$, the homoclinic orbit $\mathbf{\Gamma_{\rm HOM}}$ (brown curve) formed by one branch of $W^s(\mathbf{0})$, the other branch of $W^s(\mathbf{0})$ (cyan curve), a first part of $W^{u}(\mathbf{0})$ (red surface), and $W^{uu}(\mathbf{0})$ (magenta curve). Here $(\alpha, \beta) = (5.8, 1.7010)$ in panel $\mathbf{H_{o}}$, $(\alpha, \beta) = (5.3573, 2.1917)$ in panel $\mathbf{C}_{\rm in}$ and $(\alpha, \beta) = (5.1, 2.717)$ in panel $\mathbf{H_{t}}$.}
\end{center}
\end{figure}
%

\Fref{fig:flipTrans} illustrates the transition through the orbit flip bifurcation along the curve $\mathbf{H_{o/t}}$. On both sides of $\mathbf{C}_{\rm in}$, the one-dimensional stable manifold $W^s(\mathbf{0})$ returns to $\mathbf{0}$ tangent to the weak stable eigendirection to form the homoclinic orbit $\mathbf{\Gamma_{\rm HOM}}$. At the same time, the two-dimensional unstable manifold $W^{u}(\mathbf{0})$ returns back to the saddle $\mathbf{0}$ and closes up along the one-dimensional strong unstable manifold $W^{uu}(\mathbf{0}) \subset W^{u}(\mathbf{0})$. The shown part of the surface consists of a family of orbit segments that start at distance of $10^{-3}$ from $\mathbf{0}$; it has been computed with the boundary-value problem set-up from \cite{And1,And2018,redbook}. The two typical cases of homoclinic bifurcation are that $W^{u}(\mathbf{0})$ either forms a cylinder along $\mathbf{H_{o}}$ or a M{\"o}bius strip along  $\mathbf{H_{t}}$, depending on which side of $W^{uu}(\mathbf{0})$ the stable manifold $W^s(\mathbf{0})$ returns. This is illustrated in \fref{fig:flipTrans} by the different positions on the surface $W^{u}(\mathbf{0})$ of the curve $W^{uu}(\mathbf{0})$ relative to the homoclinic orbit; see especially the enlargements. The change in orientability occurs at the point $\mathbf{C}_{\rm in}$ when $W^s(\mathbf{0})$ returns to $\mathbf{0}$ exactly along  $W^{uu}(\mathbf{0})$, which is represented in \fref{fig:flipTrans} by the respective branches of the two manifolds coinciding in panel $\mathbf{C}_{\rm in}$. As a result, the surface $W^{u}(\mathbf{0})$ comes back tangent to the strong direction and so is neither orientable nor non-orientable.

The top-left region of the bifurcation diagram in \fref{fig:bifDiagram}(b), to the left of $\mathbf{SNP}$ and above $\mathbf{H_t}$, is the only region where system~\eqref{eq:WorkEquation} has no periodic orbits as a result of the flip bifurcation. Upon crossing $\mathbf{H_t}$, a single saddle periodic orbit $\Gamma_t$ is created, which is non-orientable; hence, it has negative nontrivial Floquet multipliers. When followed around the point $\mathbf{C}_{\rm in}$, the periodic orbit $\Gamma_t$ persists throughout the different regions in the bifurcation diagram until the curve $\mathbf{PD}$, where it merges with a repelling period-doubled orbit in a subcritical period-doubling bifurcation. This turns $\Gamma_t$ into an attracting periodic orbit, which exists in the region between the curves $\mathbf{PD}$ and $\mathbf{SNP}$. Since $\Gamma_t$ is now attracting, it can transform from a non-orientable to an orientable periodic orbit, which allows it to bifurcate at $\mathbf{SNP}$ with the orientable saddle periodic orbit $\Gamma_o$ that is created upon crossing $\mathbf{H_o}$ into the region with $\hat{\beta} > 0$. 

Many more periodic orbits are created and disappear again near the orbit flip point $\mathbf{C}_{\rm in}$, and we now turn our attention to an associated global feature of the bifurcation diagram: the nature of the curves $\mathbf{H^n}$ of $n$-homoclinic bifurcations. Observe in \fref{fig:bifDiagram}(b) that each of the curves $\mathbf{H^2}$ to $\mathbf{H^6}$ emanating from $\mathbf{C}_{\rm in}$ has a fold (a maximum) with respect to $\alpha$ and then extends towards decreasing $\alpha$ and $\hat{\beta}$, past the $\alpha$-value of the point $\mathbf{C}_{\rm in}$. Hence, all these curves also exist on the side of $\mathbf{H_t}$. The curve $\mathbf{PD^2}$ emanating from $\mathbf{C}_{\rm in}$ ends on the curve $\mathbf{H^2}$ at a codimension-two orbit flip bifurcation point $\mathbf{C^2_{\rm O}}$, quite close to the fold. We find that the bifurcation diagram in the $(\alpha, \beta)$-plane is even more complicated than was suggested in~\cite{Alg2019}. We identify codimension-two inclination flip bifurcation points $\mathbf{C^n_{\rm I}}$ on each of the curves $\mathbf{H^2}$ to $\mathbf{H^6}$, again very close to where they have a fold with respect to $\alpha$; see the enlargement \fref{fig:bifDiagram}(c). Also shown in all panels is the curve $\mathbf{SNP^3}$ of saddle-node bifurcation of periodic orbits that emanates from $\mathbf{C^3_{\rm I}}$. We observe that for sufficiently small values of $\alpha$ the $n$-homoclinic orbits along the curves $\mathbf{H^3}$ to $\mathbf{H^6}$ are non-orientable.

%
\begin{figure}[t!]
\begin{center}
  \includegraphics{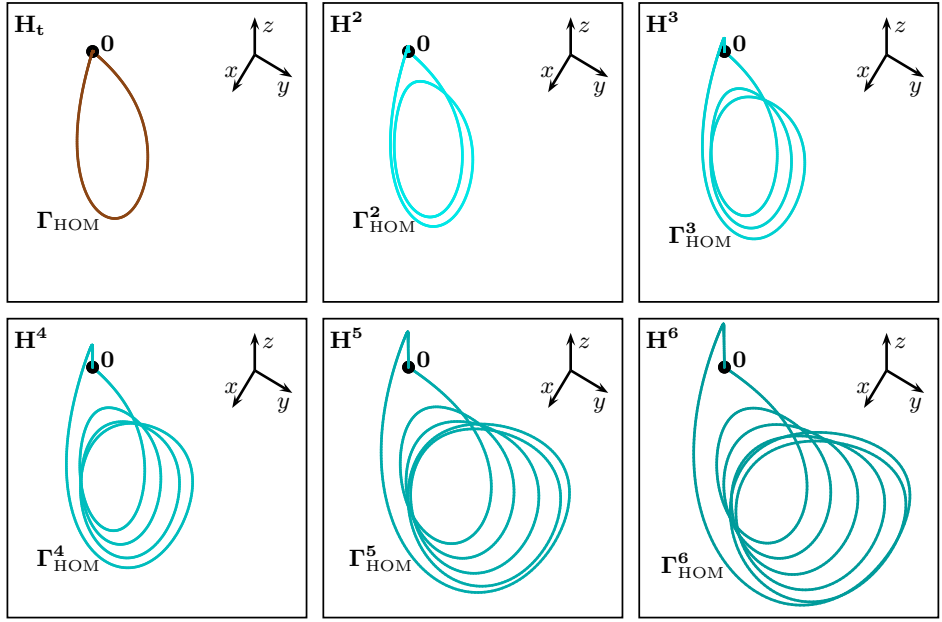}
  \caption{\label{fig:ConfigHom}
    The primary homoclinic orbit on $\mathbf{H_t}$ and the $n$-homoclinic orbits $\mathbf{H^2}$ to $\mathbf{H^6}$ of~\eqref{eq:WorkEquation} for $\alpha = 5.3$, shown in $\R^3$ in brown and increasingly darker shades of cyan to match the colors of the corresponding bifurcation curves in \fref{fig:bifDiagram}.} 
\end{center}
\end{figure}
%

The computed curves $\mathbf{H^2}$ to $\mathbf{H^6}$ in \fref{fig:bifDiagram} suggest that they are part of a family of curves $\mathbf{H^n}$ that accumulate on a well-defined limiting curve. Therefore, we now focus on the limiting behavior of the curves $\mathbf{H^n}$ and of the associated $n$-homoclinic orbits as the number of loops $n$ increases. The homoclinic orbits on $\mathbf{H_t}$ and $\mathbf{H^2}$ to $\mathbf{H^6}$ for fixed $\alpha = 5.3$ are shown in \fref{fig:ConfigHom}, where they are assigned the same color as the corresponding curves in \fref{fig:bifDiagram}. Each time, from one panel to the next, the branch of $W^s(\mathbf{0})$ that forms the homoclinc orbit has one extra loop before closing up. Notice that, with increasing $n$, the additional loops of the homoclinic orbit extend increasingly futher along the $y$-direction. This behavior is intriguing, because it suggest that the $n$-homoclinic orbits converge with $n$ to a heteroclinic connection from $\mathbf{0}$ to an equilibrium or periodic orbit at infinity, which corresponds to the limiting case of an infinite number of larger and larger loops. This suggests that the curves $\mathbf{H^n}$ in the two-parameter plane accumulate onto a curve of such heteroclinic bifurcations involving infinity, which is, therefore, expected to be of codimension one.

\section{Characterizing the dynamics at infinity}
\label{sec:compact}

For the purpose of finding a possible heteroclinic bifurcation involving infinity, we must identify equilibria or periodic orbits at infinity. We take advantage of the fact that system~\eqref{eq:WorkEquation} is a polynomial vector field, which means that we can compactify the phase space. In general terms, the behavior at infinity is given, after a suitable compactification, by the terms of highest order. We identify and analyze different invariant objects in new coordinate charts that represent the dynamics at and near infinity. This approach makes it possible to continue equilibria or other special solutions as they interact in degenerate bifurcations at infinity~\cite{And1}. The purpose here is to use charts at infinity to set up a well-posed boundary value problem with a solution that represents the heteroclinic connection to infinity. 

More specifically, we follow the recent work by Matsue~\cite{Matsue2018} to obtain a suitable Poincar{\'e} compactification for system~\eqref{eq:WorkEquation}; see also~\cite{Vel, Jes1, Jes2}. The underlying idea was already proposed in~\cite{Dumortier1991} for planar vector fields, where it is defined as a directional blow-up for so-called quasi-homogeneous vector fields. In our context, this means applying a directional compactification in the direction of positive $y$, because the $n$-homoclinic orbits extend predominantly in the $y$-direction as $n$ increases, while their $x$- and $z$-components remain relatively bounded. 

Note that system~\eqref{eq:WorkEquation}  is not quasi-homogeneous. However, investigation of the leading terms of the right-hand side of system~\eqref{eq:WorkEquation} in the limit to infinity shows that it is \emph{asymptotically quasi-homogeneous}~\cite{Matsue2018} to the quasi-homogeneous vector field of type $(3, 4, 1)$ and order $3$ given by
\begin{equation}
\label{eq:quasihom}
  \left\{ \begin{array}{rcl}
            \dot{x} &=& y \, z, \\
            \dot{y} &=& x^2, \\
            \dot{z} &=& -4 x.
          \end{array} \right.
\end{equation} 
The powers of the directional blow-up are then determined by the type of the quasi-homogeneous system~\eqref{eq:quasihom}, which leads to the coordinate transformation
\begin{displaymath}
  (x, y, z) \mapsto (\bar{x}, \bar{z}, \bar{w}), \quad
  x = \frac{\bar{x}}{\bar{w}^3},\; y = \frac{1}{\bar{w}^4}, \mbox{ and } z = \frac{\bar{z}}{\bar{w}}.
\end{displaymath}
These coordinates define the chart with $y > 0$, and $\bar{w}$ represents the distance to infinity in the $y$-direction. More precisely, let $(x_{\rm s},y_{\rm s},z_{\rm s})$ be the transformed coordinates of system~\eqref{eq:quasihom} inside the Poincar{\'e} sphere centered at the origin, where directions of escape to infinity are represented  by points on the sphere of radius one. In these coordinates, $(\bar{x},\bar{y},\bar{z})$ correspond to the projection of the positive $y_{\rm s}$-hemisphere of the two-dimensional Poincar{\'e} sphere onto the plane defined by $y_{\rm s} = 1$. The resulting weighted directional compactification then becomes
\begin{displaymath}
  \left\{ \begin{array}{rcrl}
            \dot{\bar{x}} &=& \dfrac{1}{\bar{w}^2} & \left( \bar{z} + \alpha \, \bar{w} +
                  \frac{3}{4} \bar{x} \, [\bar{w}^2 - \beta \, \bar{w}^3 \, \bar{x} - \bar{x}^2] + \gamma \, \bar{z} \, \bar{w}^4 \right), \\[3mm]
            \dot{\bar{z}} &=& \dfrac{1}{\bar{w}^2} & \left( -4 \bar{x} +
                  \frac{1}{4} \bar{z} \, [\bar{w}^2 - \beta \, \bar{w}^3 \, \bar{x} - \bar{x}^2] \right), \\[3mm]
            \dot{\bar{w}} &=& \dfrac{1}{\bar{w}^2} & \left( \frac{1}{4} \bar{w} \, [\bar{w}^2 - \beta \, \bar{w}^3 \, \bar{x} - \bar{x}^2] \right).
          \end{array} \right.
\end{displaymath} 
It can be desingularized via a rescaling of time with the factor $\bar{w}^2$, yielding the desingularized vector field that contains 
the dynamics at infinity as
\begin{equation}
\label{eq:ChartEquation}
  \left\{ \begin{array}{rcl}
            \dot{\bar{x}} &=& \bar{z} + \alpha \, \bar{w} +
              \frac{3}{4} \bar{x} \,  (\bar{w}^2 - \beta \, \bar{w}^3 \, \bar{x} - \bar{x}^2) + \gamma \, \bar{z} \, \bar{w}^4, \\[2mm]
            \dot{\bar{z}} &=& -4 \bar{x} +
                  \frac{1}{4} \bar{z} \, (\bar{w}^2 - \beta \, \bar{w}^3 \, \bar{x} - \bar{x}^2), \\[2mm]
            \dot{\bar{w}} &=& \frac{1}{4} \bar{w} \, (\bar{w}^2 - \beta \, \bar{w}^3 \, \bar{x} - \bar{x}^2).
      \end{array} \right.
\end{equation} 
\begin{remark}
It is also possible to perform a standard directional Poincar\'{e} compactification that gives all variables the same weight. However, we found that this leads to highly non-hyperbolic dynamics in the chart with $y > 0$ so that the dynamics at infinity is difficult to characterize. This issue would then have to be resolved via a blow-up procedure  with exponents that take into account the weighting used to obtain system~\eqref{eq:ChartEquation}.
\end{remark}

%
\begin{figure}[t!]
\begin{center}
  \includegraphics{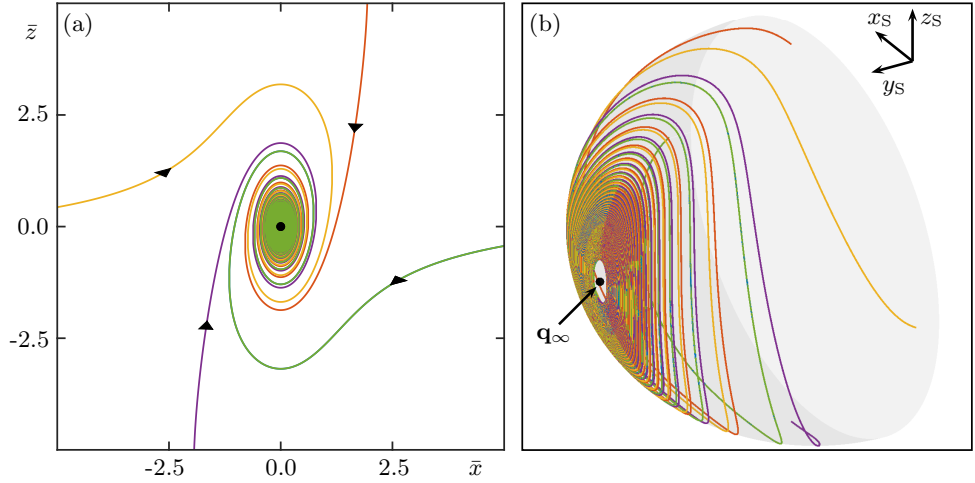}
  \caption{\label{fig:Poincare}
    Dynamics at infinity for system~\eqref{eq:ChartEquationZero}, or system~\eqref{eq:ChartEquation} with $\bar{w}=0$, shown in the $(\bar{x}, \bar{z})$-plane in panel~(a). Panel~(b) shows the projection of panel~(a) onto the corresponding Poincar{\'e} half-sphere with $y_{\rm s} > 0$ in the compactified $(x_{\rm s},y_{\rm s},z_{\rm s})$-cordinates.}
  \end{center}
\end{figure}
%

We are now ready to analyze the dynamics at infinity and decide whether it contains equilibria or periodic orbits that could be involved in a suspected heteroclinic connection. To this end, we set $\bar{w}=0$ in system~\eqref{eq:ChartEquation} and observe that the $(\bar{x}, \bar{z})$-plane is indeed invariant. The resulting system
\begin{equation}
\label{eq:ChartEquationZero}
  \left\{ \begin{array}{rcrcl}
            \dot{\bar{x}} &=& \bar{z} &-& \frac{3}{4} \bar{x}^3, \\[2mm]
            \dot{\bar{z}} &=& -4 \bar{x} &-& \frac{1}{4} \bar{x}^2 \, \bar{z},
          \end{array} \right.
\end{equation}
has a single equilibrium at $(\bar{x}, \bar{z}) = (0, 0)$, which is, in fact, not hyperbolic. This equilibrium is the equilibrium $\mathbf{q}_\infty$ at infinity in system~\eqref{eq:WorkEquation}. To understand the dynamics at infinity, that is, on the $(\bar{x}, \bar{z})$-plane, we convert to polar coordinates. More precisely, we consider the ellipsoidal transformation
\begin{displaymath}
  (\bar{x}, \bar{z}) \mapsto (\bar{r},\bar{\theta}), \quad
  \bar{x} = \bar{r} \, \cos{\bar{\theta}}, \mbox{ and }\bar{z} = 2 \bar{r} \, \sin{\bar{\theta}}.
\end{displaymath}
The vector field in these ellipsoidal polar coordinates becomes 
\begin{displaymath} 
  \left\{ \begin{array}{rcl}
            \dot{\bar{r}} &=&  -\frac{1}{4}\bar{r}^3\cos^2{\bar{\theta}} \left( 2+\cos{2 \bar{\theta}} \right), \\[2mm]
            \dot{\bar{\theta}} &=& -2+\tfrac{1}{2}\bar{r}^2\cos^3{\bar{\theta}}\sin{\bar{\theta}}.
      \end{array} \right.
\end{displaymath} 
Note that $\dot{\bar{r}} < 0$ for all $(\bar{r}, \bar{\theta})$ with $\bar{r} > 0$, and that $\dot{\bar{\theta}} < 0$ and close to $-2$ as soon as $\bar{r}$ is small enough. Hence, all trajectories in the $(\bar{x}, \bar{z})$-plane converge to $\mathbf{q}_\infty$, which lies at the origin in this planar coordinate system; moreover, locally near $\mathbf{q}_\infty$, trajectories will spiral clockwise towards it. This behavior is illustrated in \fref{fig:Poincare}, where we plot several trajectories in the $(\bar{x}, \bar{z})$-plane in panel~(a) and project them back onto the Poincar{\'e} sphere in panel~(b); note that system~\eqref{eq:ChartEquationZero} only describes the dynamics in the chart with $y_{\rm s} > 0$, and only the corresponding half-sphere is shown. 

%
\begin{figure}[t!]
\begin{center}
  \includegraphics{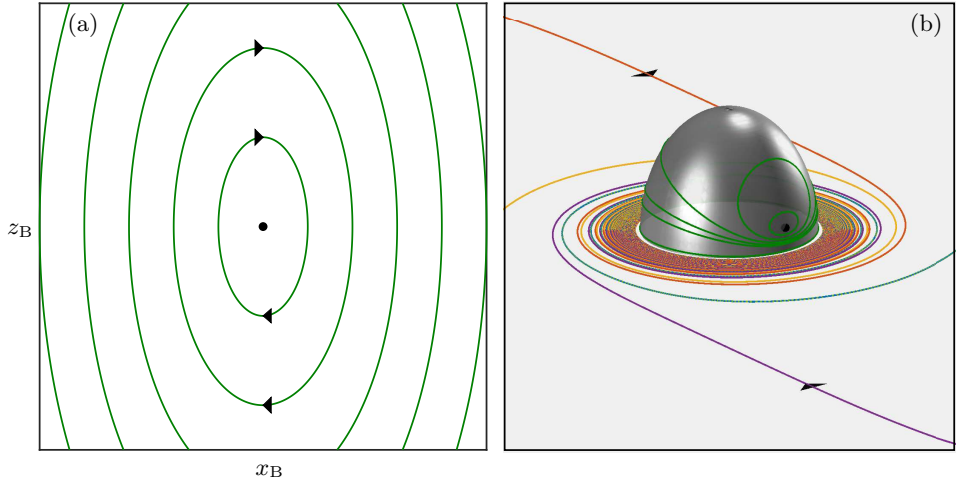}
  \caption{\label{fig:Blow}  
    Dynamics near the equilibrium $(\bar{x}, \bar{z}, \bar{w}) = (0, 0, 0)$ of system~\eqref{eq:ChartEquation}. The behavior in the $(x_{\rm B}, z_{\rm B})$-plane, that is, the blow-up chart \eqref{eq:BlowUpEquation} with $w_{\rm B} = 0$, is shown in panel~(a). It corresponds to the dynamics on a half-sphere around the origin in the $(\bar{x}, \bar{z}, \bar{w})$-space, as is illustrated in panel~(b); compare also with \fref{fig:Poincare}(a).}    
  \end{center}
\end{figure}
%

In the full three-dimensional blown-up system~\eqref{eq:ChartEquation}, the point $\mathbf{q}_\infty$ at $(\bar{x}, \bar{z}, \bar{w}) = (0, 0, 0)$ is not a hyperbolic attractor. Therefore, we perform an additional $\bar{w}$-directional blow-up by applying the transformation 
\begin{displaymath}
  (\bar{x}, \bar{z}, \bar{w}) \mapsto (x_{\rm B},z_{\rm B},w_{\rm B}), \quad
  \bar{x} = x_{\rm B} \, w_{\rm  B},\; \bar{z} = z_{\rm B} \, w_{\rm  B}, \mbox{ and } \bar{w} = w_{\rm  B},
\end{displaymath}
to system~\eqref{eq:ChartEquation}. This gives the vector field
\begin{equation}
\label{eq:BlowUpEquation}
  \left\{ \begin{array}{rcl}
       \dot{x}_{\rm B} &=& \alpha + z_{\rm B} + \gamma \, w_{\rm B}^4 \, z_{\rm B} 
                           + \frac{1}{2} x_{\rm B} \, w_{\rm B}^2 \, (1 - \beta \, x_{\rm B} \, w_{\rm  B}^2 - x_{\rm B}^2), \\[2mm]
       \dot{z}_{\rm B}  &=& -4 x_{\rm B},\\[2mm]
       \dot{w}_{\rm B}  &=& \frac{1}{4} w_{\rm B}^3 \, (1 - \beta x_{\rm B} \, w_{\rm B}^2 - x_{\rm B}^2),
   \end{array} \right.
\end{equation} 
which further characterizes the dynamics at infinity on a local half-sphere around $\mathbf{q}_\infty$. Setting $w_{\rm B}=0$ in system~\eqref{eq:BlowUpEquation}, we find that the invariant $(x_{\rm B}, z_{\rm B})$-plane is foliated by ellipses of the form $4 x^2_{\rm B} + (z_{\rm B} + \alpha)^2 = c^2$; see \fref{fig:Blow}(a). The trajectories in the $(x_{\rm B}, z_{\rm B})$-plane correspond to trajectories on the blown-up half-sphere with $\bar{w} > 0$ centered at $(\bar{x}, \bar{z}, \bar{w}) = (0, 0, 0)$. \Fref{fig:Blow}(b) gives an impression of how the previously identified dynamics at infinity interacts with the blown-up half-sphere in the $(x_{\rm B}, z_{\rm B}, w_{\rm B})$-space.

\begin{figure}[t!]
\begin{center}
  \includegraphics{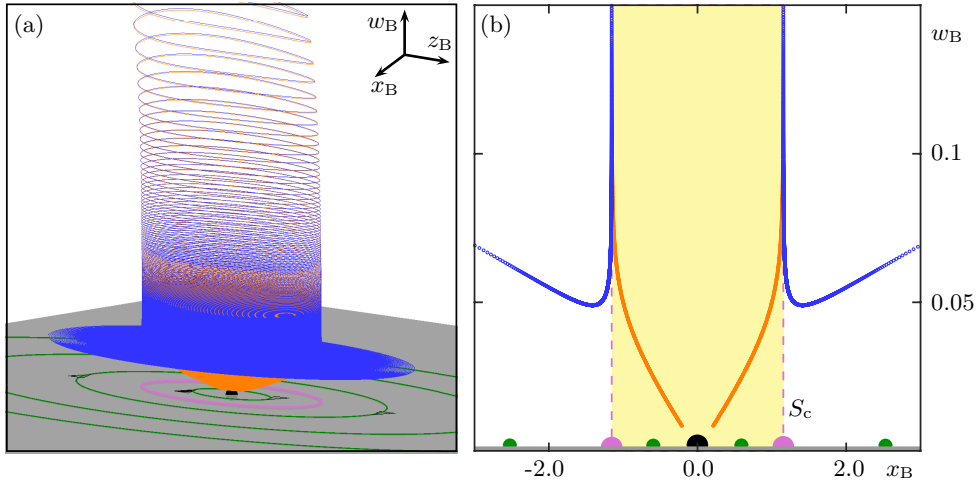}
  \caption{\label{fig:Simulations}  
    Numerical simulations suggest the existence of a cylinder-shaped separatrix $S_{\rm c}$ of system~\eqref{eq:BlowUpEquation} between trajectories that converge to the equilibrium $(x_{\rm B}, z_{\rm  B}, w_{\rm B}) = (0, -\alpha, 0)$, such as the orange trajectory, and those that do not, such as the blue trajectory. Panel~(a) shows the $(x_{\rm B}, z_{\rm B}, w_{\rm B})$-space near $(0, -\alpha, 0)$ and panel~(b) the associated intersection sets with with the plane defined by $z_{\rm B}= -\alpha$.}
  \end{center}
\end{figure}

The next step is to determine the property of system~\eqref{eq:BlowUpEquation} for $\bar{w} > 0$. First, we resort to numerical simulation and determine how initial conditions with $\bar{w} > 0$ approach the $(x_{\rm B}, z_{\rm B})$-plane. \Fref{fig:Simulations} shows that there are two types of behavior. Panel~(a) shows two trajectories of system~\eqref{eq:BlowUpEquation} for $(\alpha,\beta)=(5.3,2.0)$, obtained by integration in both forward and backward time from the initial conditions $(x_{\rm B}, z_{\rm  B}, w_{\rm B}) = (1, -\alpha, 0.05)$ and $(x_{\rm B}, z_{\rm  B}, w_{\rm B}) = (1.3, -\alpha, 0.05)$, respectively. The former initial condition leads to a trajectory (orange) that converges in backward time in a spiraling fashion to the equilibrium $(x_{\rm B}, z_{\rm  B}, w_{\rm B}) = (0, -\alpha, 0)$ of~\eqref{eq:BlowUpEquation}. The other trajectory (blue) first approaches the $(x_{\rm B}, z_{\rm B})$-plane in backward time but then diverges away from it; in particular, it does not reach the equilibrium $(x_{\rm B}, z_{\rm  B}, w_{\rm B}) = (0, -\alpha, 0)$. This is illustrated further in \fref{fig:Simulations}(b) in a local cross-section defined by $z_{\rm B} = -\alpha$. Notice that the two trajectories are very close together before they separate in backward time at about $\bar{w} = 0.08$. 

%
%

We conclude that there exists an invariant critical surface $S_{\rm c}$ that separates the two qualitatively different regions in phase space where trajectories converge to $(x_{\rm B}, z_{\rm  B}, w_{\rm B}) = (0, -\alpha, 0)$ and where they do not. Furthermore, \fref{fig:Simulations} suggests that this difference in the backward-time limit of trajectories is entirely due to the fact that ellipses near $(x_{\rm B}, z_{\rm  B}, w_{\rm B}) = (0, -\alpha, 0)$ are repelling in the $\bar{w}$-direction, while beyond some distance, they are attracting in the $\bar{w}$-direction. The surface $S_{\rm c}$ is associated with the critical ellipse in the $(x_{\rm B}, z_{\rm B})$-plane that is neither repelling nor attracting in the $\bar{w}$-direction, and which goes through the point $(x_{\rm B}, z_{\rm   B}) \approx (1.1547, -\alpha)$ (magenta). Our computations indicate that the critical surface $S_{\rm c}$ is effectively a straight elliptical cylinder when $\bar{w}$ is small. 

Based on these careful observations, we approximate $S_{\rm c}$ as a straight $\bar{w}$-cylinder around the $w_B$-axis through the point $(0, -\alpha, 0)$.  We denote this cylinder $C_{r^*}$, with a specific radius $r^*$, and require that
\begin{equation}
\label{eq:zeroflux}
  \int_{C_{r}} {X_{\rm B} \cdot \hat{n}_{C_{r}} \; d C_{r}} = 0
\end{equation}
be satisfied for $r=r^*$, where $\hat{n}_{C_{r}}$ is the direction normal to $C_{r}$. We are using the average zero-flux condition  \eqref{eq:zeroflux} to define $C_{r^*}$ because, in general, there is no cylinder that is invariant under the vector field $X_{\rm B}$ defined by~\eqref{eq:BlowUpEquation}. To find  $r^*$ we transform system~\eqref{eq:BlowUpEquation} to cylindrical coordinates by
\begin{displaymath}
  (x_{\rm B}, z_{\rm B}, w_{\rm  B}) \mapsto (r_{\rm B}, \theta_{\rm B}, w_{\rm B}), \quad 
  x_{\rm B} = r_{\rm B} \, \cos{\theta_{\rm  B}} \mbox{ and } z_{\rm B} = 2 r_{\rm B} \, \sin{\theta_{\rm  B}} - \alpha.
\end{displaymath}
The integral can then be evaluated in a straightforward way as:
\begin{eqnarray*}
  \int_{C_r} {X_{\rm B} \cdot \hat{n}_{C_r} \; d C_r} 
  &=& \left. \int d w_{\rm B} \; \int_0^{2\pi} \dot{r}_{\rm B} \; d \theta_{\rm B} \right|_{r_{\rm B} = r} \\ 
  &=& \int \frac{1}{8} \pi r \, w_{\rm B}^2 \, (4 - 3 r^2) \; d w_{\rm B} 
  = \frac{1}{24} \pi r \, w_{\rm B}^3 \, (4 - 3 r^2).
\end{eqnarray*}
%
Hence, there are two zeros of the zero-flux condition~\eqref{eq:zeroflux}, namely, $r= 0$ and $r=r^* = \frac{2}{3} \sqrt{3}$. Note that $r = 0$ corresponds to the $w_{\rm B}$-axis through the equilibrium at $(0, -\alpha, 0)$. We conclude that the critical cylinder $C_{r^*}$ with $r^* = \frac{2}{3} \sqrt{3} \approx 1.1547$ is the local approximation of the separating invariant surface $S_{\rm c}$. This value agrees with our numerical simulations and $C_{r^*}$ is a good first-order approximation of $S_{\rm c}$.

%
\begin{figure}[t!]
\begin{center}
  \includegraphics{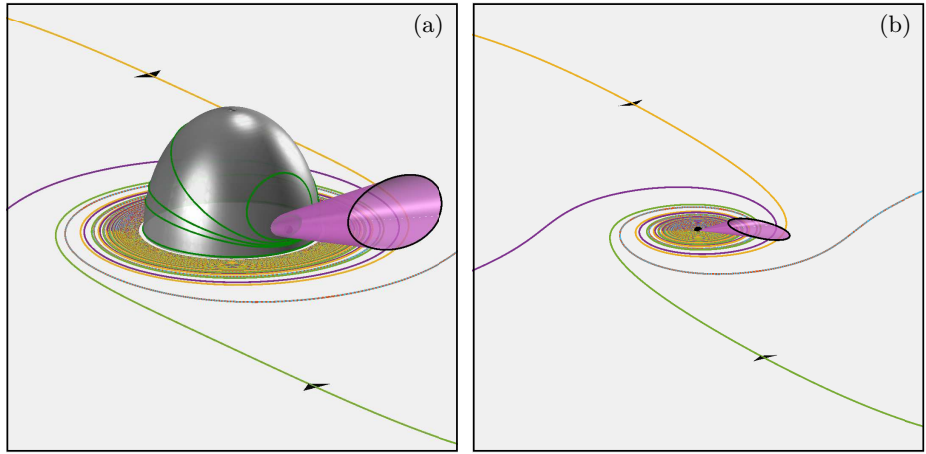}
  \caption{\label{fig:deBlowup}
    The separatrix $S_{\rm c}$ (purple surface) as represented locally by the cylinder $C_{r^*}$, shown in the $(\bar{x}, \bar{z}, \bar{w})$-space of system~\eqref{eq:ChartEquation}. Panel~(a) shows $S_{\rm c}$ emerging from the blown-up half-sphere, while in panel~(b), $S_{\rm c}$ is a cone that emerges from the origin.}
  \end{center}
\end{figure}
%

Recall that the $(x_{\rm B}, z_{\rm B}, w_{\rm B})$-coordinate system of~\eqref{eq:BlowUpEquation} corresponds to a directional blow-up of the equilibrium $\mathbf{q}_\infty$ at infinity in the original coordinates, which corresponds to the origin in the $(\bar{x}, \bar{z}, \bar{w})$-coordinates of the desingularized system~\eqref{eq:ChartEquation}. \Fref{fig:deBlowup} illustrates in two ways the separatrix $S_{\rm c}$ (magenta surface) represented by the inverse image of the critical cylinder $C_{r^*}$ under the respective coordinate transformations. Panel~(a) shows how $S_{\rm c}$ emanates from a corresponding periodic orbit on the blown-up half-sphere centered at the origin of the $(\bar{x}, \bar{z}, \bar{w})$-space. However, periodic orbits only exist on the blown-up (half-)sphere and not in the $(\bar{x}, \bar{z}, \bar{w})$-space itself. Deflating the blown-up sphere back to the origin, the local approximation $C_{r^*}$  of $S_{\rm c}$ is the cone emanating from the origin in the $(\bar{x}, \bar{z}, \bar{w})$-space that is shown in \fref{fig:deBlowup}(b).

\section{BVP set-up for computing a codimension-one connection to infinity}
\label{sec:LinAtInfty}
All trajectories inside the separatrix $S_{\rm c}$ converge, in backward time, to $\mathbf{q}_\infty$, which is the origin in the $(\bar{x}, \bar{z}, \bar{w})$-space. Hence, $S_{\rm c}$ acts as a kind of two-dimensional unstable manifold of the non-hyperbolic point $\mathbf{q}_\infty$ at infinity. For special choices of the parameters $\alpha$ and $\beta$ in system~\eqref{eq:WorkEquation}, the one-dimensional stable manifold $W^s(\mathbf{0})$ of the origin in the original $(x, y, z)$-coordinates lies in the surface $S_{\rm c}$. We refer to this well-defined phenomenon of codimension one as a heteroclinic connection between $\mathbf{0}$ and $\mathbf{q}_\infty$, and we denote it by $\mathbf{Het^\infty}$. It is our hypothesis that the curves $\mathbf{H^n}$ of $n$-homoclinic orbits, which have increasingly longer excursions towards infinity, accumulate in the $(\alpha,\beta)$-plane on the corresponding curve $\mathbf{Het^\infty}$; see \fref{fig:ConfigHom}. 

Hence, the task is to find the heteroclinic connection $\mathbf{Het^\infty}$ and to continue it in the $(\alpha,\beta)$-plane. To this end, we employ the approach known as Lin's method~\cite{KraRie1, Lin1990} to set up a two-point boundary value problem (BVP) for two orbit segments such that their concatenation is the sought-after connecting orbit in $W^s(\mathbf{0}) \cap S_{\rm c}$. The essence of Lin's method is to choose a codimension-one plane $\Sigma$ that separates the two invariant objects involved, here $\mathbf{0}$ and $\mathbf{q}_\infty$, and to consider an orbit segment in $W^s(\mathbf{0})$ up to $\Sigma$ and an orbit segment in $S_{\rm c}$ up to $\Sigma$. For parameters that are not at the bifurcation value, these two orbit segments exhibit a gap in $\Sigma$. Lin's theorem states that this orbit pair and, hence, the gap are uniquely determined when the difference between their end points in $\Sigma$ is constrained to lie in a fixed subspace called the Lin space~\cite{Lin1990}. The associated signed Lin gap in the Lin space is then a well-defined test function with zeros that correspond to connecting orbits; such zeros can be found via the continuation of the corresponding orbit segements as solutions of an overall BVP~\cite{KraRie1}. Once a zero is found, the associated connecting orbit can be followed in system parameters.

The challenge, here, is that one of the equilibria lies at infinity and we have an approximation for $S_{\rm c}$ in blown-up coordinates. Note that systems~\eqref{eq:WorkEquation} and~\eqref{eq:BlowUpEquation} are homeomorphic in the open sets where they coincide~\cite{Matsue2018}. This allows us to define $\Sigma$ with respect to both coordinate systems. We then consider one orbit segment that is a solution of system~\eqref{eq:WorkEquation} with one end point near the saddle $\mathbf{0}$ and lying in its stable eigenspace (which is the linear approximation of $W^s(\mathbf{0})$) and the other lying in $\Sigma$; and a second orbit segment that is a solution of system~\eqref{eq:BlowUpEquation} with one end point near the point $(\bar{x},\bar{y}, \bar{z}) = (0,0,0)$ representing $\mathbf{q}_{\infty}$  lying in the linear approximation $C_r$ of $S_{\rm c}$ and the other lying in $\Sigma$. The respective coordinate transformations allow us to `glue' the original $(x, y, z)$-coordinates of  system~\eqref{eq:WorkEquation} to the $(x_{\rm B}, z_{\rm B}, w_{\rm B})$-coordinates of the blown-up system~\eqref{eq:BlowUpEquation}, so that we can define and determine the Lin gap.

We use this adapted Lin's method to find an initial connecting orbit in $W^s(\mathbf{0}) \cap S_{\rm c}$, along with the relevant bifurcation value for $\beta$, where we keep $\alpha = 5.3$ fixed. We define 
\begin{displaymath}
  \Sigma := \{ (x,y,z) \;|\; x = 0 \} \simeq \{ (x_{\rm B}, z_{\rm B}, w_{\rm B}) \;|\; x_{\rm B} = 0 \},
\end{displaymath}
which is a suitable choice that works in both coordinate systems for $\bar{w} \neq 0$ because $x = \bar{x} / \bar{w}^3$ and $\bar{x} = x_{\rm B} \, w_{\rm  B}$ with $w_{\rm  B} = \bar{w}$. 

To define the orbit segment $\mathbf{u}$ in $(x,y,z)$-coordinates that lies in $W^s(\mathbf{0})$ up to $\Sigma$ we define the BVP
\begin{eqnarray}
\label{eq:BVPto0_vf}
            \dot{\mathbf{u}} &=& T_0 \, X(\mathbf{u}), \\
\label{eq:BVPto0_zero}
            \mathbf{u}(1) &=& \delta_0 \, \mathbf{e}^s, \\
\label{eq:BVPto0_sigma}
            \mathbf{u}(0)^* \, \mathbf{n}  &=& 0.
\end{eqnarray}
Here, $X$ denotes the vector field~\eqref{eq:WorkEquation} and $T_0$ is the total integration time between the first and last point on the orbit segment; it enters \eqref{eq:BVPto0_vf} in explicit form so that the orbit segment $\mathbf{u}(t)$ is defined for $t \in [0,1]$. Boundary condition~\eqref{eq:BVPto0_zero} requires that the end point $\mathbf{u}(1)$ lies at a small distance $\delta_0$ from the saddle $\mathbf{0}$ along its stable eigenvector $\mathbf{e}^s$ (which has been normalized to have length 1). This ensures that $\mathbf{u}(1)$ lies in $W^s(\mathbf{0})$ to good approximation, provided $\delta_0$ is sufficiently small; we fix $\delta_0 = 10^{-4}$ as an appropriate value throughout. Finally, the dot product in boundary condition~\eqref{eq:BVPto0_sigma} involves the unit vector $\mathbf{n} = (1, 0, 0)$ normal to $\Sigma$, which ensures that the start point $\mathbf{u}(0)$ lies in $\Sigma$. We remark that the stable eigenvector $\mathbf{e}^s$ in ~\eqref{eq:BVPto0_zero} needs to be continued as well when system parameters are changed; we achieve this by solving the BVP of the corresponding stable eigenvector problem~\cite{KraRie1} together with~\eqref{eq:BVPto0_vf}--\eqref{eq:BVPto0_sigma}. 

Similarly, the orbit segment $\mathbf{u}$ in $(x_{\rm B}, z_{\rm B}, w_{\rm B})$-coordinates in $S_{\rm c}$ up to $\Sigma$ is defined by the BVP
\begin{eqnarray}
\label{eq:VFtoinfty_vf}
  \dot{\mathbf{u}}_{\rm B} &=& T_{\rm B} \, X_{\rm B}(\mathbf{u}_{\rm B}), \\
\label{eq:BVPtoinfty_zero}
  \mathbf{u}_{\rm B}(0) &=& ( \tfrac{2}{3} \sqrt{3} \, \cos{\theta_{\rm B}},\, \tfrac{4}{3} \sqrt{3} \, \sin{\theta_{\rm B}} - \alpha,\, \delta_{\rm B} ), \\
\label{eq:BVPtoinfty_sigma}
            \mathbf{u}_{\rm B}(1)^* \, \mathbf{n}  &=& 0.
\end{eqnarray}
In \eqref{eq:VFtoinfty_vf} the vector field~\eqref{eq:BlowUpEquation} is denoted $X_{\rm B}$, and $T_{\rm B}$ is the total integration time. Boundary condition \eqref{eq:BVPtoinfty_zero} requires that the start point $\mathbf{u}_{\rm B}(0)$ lies in the cylinder $C_{r^*}$, which has been parameterized by the angle $\theta_{\rm B} \in [0,\, 2 \pi]$ and the distance $\delta_{\rm B} $ in the $w_{\rm B}$-direction; we set $\delta_{\rm B} = 0.1$ throughout. Boundary condition~\eqref{eq:BVPtoinfty_sigma} again ensures that the end point $\mathbf{u}_{\rm B}(1)$ lies in $\Sigma$, because $\mathbf{n} = (1, 0, 0)$ is also the unit normal to $\Sigma$ in $(x_{\rm B}, z_{\rm B}, w_{\rm B})$-coordinates. 

To find first orbit segments $\mathbf{u}$ and $\mathbf{u}_{\rm B}$ that satisfy \eqref{eq:BVPto0_vf}--\eqref{eq:BVPto0_sigma} and \eqref{eq:VFtoinfty_vf}--\eqref{eq:BVPtoinfty_sigma}, respectively, we fix $\beta = 1.8$ and proceed as follows (recall that $\alpha = 5.3$ and $\gamma =0.5$ are fixed). For $\mathbf{u}$ we require initially only \eqref{eq:BVPto0_vf} and~\eqref{eq:BVPto0_zero}, and start a continuation in the integration time $T_0$ from $T_0 = 0$; note that this constitutes solving the initial value problem from the point $\delta_0 \, \mathbf{e}^s$ by continuation. During this computation we monitor the dot product and record whenever $\mathbf{u}$ satisfies condition~\eqref{eq:BVPto0_sigma}, that is, $\mathbf{u}(0)$ lies in $\Sigma$. Similarly, for $\mathbf{u}_{\rm B}$ we require only \eqref{eq:VFtoinfty_vf} and~\eqref{eq:BVPtoinfty_zero}; we start with $\theta_{\rm B} = 0$ and continue in $T_{\rm B}$ from $T_{\rm B} = 0$, while recording whenever \eqref{eq:BVPtoinfty_sigma} is satisfied and $\mathbf{u}_{\rm B}(1)$ lies in $\Sigma$. We remark that both conditions~\eqref{eq:BVPto0_sigma} and~\eqref{eq:BVPtoinfty_sigma} are satisfied for many values of $T_0$ and $T_B$, respectively, because the trajectories that contain $\mathbf{u}$ and $\mathbf{u}_B$ intersect $\Sigma$ many times.

We choose orbit segments $\mathbf{u}$ and $\mathbf{u}_{\rm B}$ that have end points in $\Sigma$ which lie suitably close to each other and couple them by defining the Lin space and associated Lin gap. To this end, we define and then fix the unit vector
\begin{displaymath}
  \mathbf{\Psi} := \frac{\mathbf{u}(0) - \widetilde{\mathbf{u}}_{\rm B}(1)}{\mid\!\mid\! \mathbf{u}(0) - \widetilde{\mathbf{u}}_{\rm B}(1) \!\mid\!\mid}
\end{displaymath}
given by the initial chosen end points of $\mathbf{u}$ and $\mathbf{u}_{\rm B}$; here $\widetilde{\mathbf{u}}_{\rm B}(1)$ is the end point of $\mathbf{u}_{\rm B}(1)$ in the original $(x,y,z)$-ccordinates of the section $\Sigma$. The vector $\mathbf{\Psi}$ is generically transverse to $S_{\rm c} \cap \Sigma$, spans the Lin space $Z$, and defines the Lin gap $\eta$ via the boundary condition
\begin{equation}
\label{eq:BVPtoinftyL}
  \widetilde{\mathbf{u}}_{\rm B}(1) = \mathbf{u}(0) + \eta \, \mathbf{\Psi}.
\end{equation}
Note that the new parameter $\eta$ is the signed distance between the two end points of the orbit segments along the Lin space $Z \subset \Sigma$, which is fixed once chosen in this way. 

We now consider the combined boundary value problem given by  \eqref{eq:BVPto0_vf}--\eqref{eq:BVPtoinfty_zero} and \eqref{eq:BVPtoinftyL}, which is automatically satisfied by the chosen orbit segments $\mathbf{u}$ and $\mathbf{u}_{\rm B}$ and uniquely defines the Lin gap $\eta$. When $\mathbf{u}$ and $\mathbf{u}_{\rm B}$ are continued in $\beta$, where $\theta_{\rm B} \in [0,\, 2 \pi]$, $T_0 > 0$, $T_{\rm B} > 0$, and $\eta \in \R$ are free parameters (but crucially $Z \subset \Sigma$ remains fixed), the Lin gap $\eta$ is monitored. When $\beta$ changes, the orbit segment $\mathbf{u}$ as well as the $\theta$-dependent orbit segment $\mathbf{u}_{\rm B}$ vary. In light of the Lin condition \eqref{eq:BVPtoinftyL}, the angle parameter $\theta_{\rm B}$ is adjusted automatically in such a way that the end point $\mathbf{u}_{\rm B}(1)$ only varies along the direction $\mathbf{\Psi}$, either away from or towards $\mathbf{u}(0)$. When a zero of $\eta$ is detected then we have found the value of $\beta$ at which the heteroclinic connection $\mathbf{Het^\infty}$ occurs; the corresponding heteroclinic orbit that connects $\mathbf{q}_\infty$ with $\mathbf{0}$ is given as the concatenation of $\mathbf{u}$ and $\mathbf{u}_{\rm B}$.

%
\begin{figure}[t!]
\begin{center}
  \includegraphics{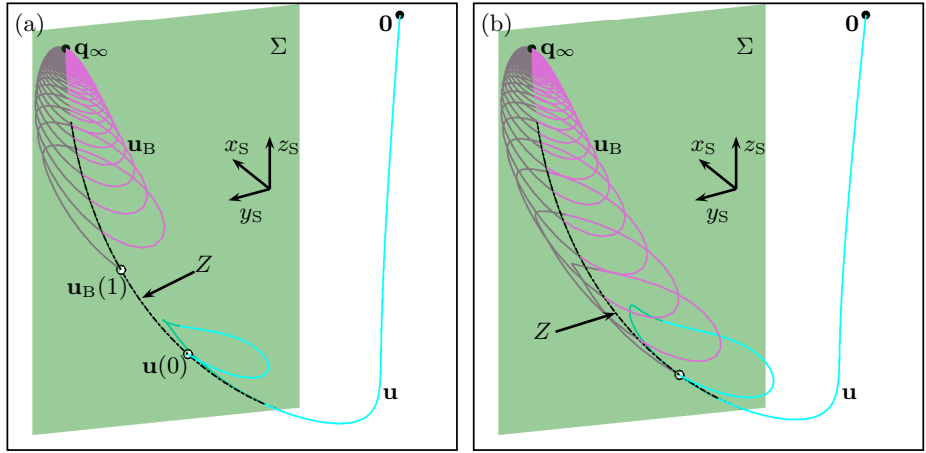}
  \caption{\label{fig:Lin}
    Set-up with Lin's method to compute a connecting orbit from $\mathbf{q}_\infty$ to $\mathbf{0}$ with two orbit segments that meet in the common Lin section $\Sigma$ (green plane), illustrated in compactified Poincar{\'e} coordinates. Panel~(a) shows the initially chosen orbit segments $\mathbf{u}$ (cyan) to $\mathbf{0}$ and $\mathbf{u}_{\rm B}$ (magenta) from $\mathbf{q}_\infty$ for $\beta = 1.8$ that define the Lin space $Z$ (which appears curved in this representation); note that the Lin gap $\eta$ is initially nonzero. Panel~(b) shows the situation for $\beta=2.08874$ where $\eta = 0$ and $\mathbf{u}$ and $\mathbf{u}_{\rm B}$ connect in $\Sigma$ to form the heteroclinic connection; here, $\alpha = 5.3$.}
\end{center}
\end{figure}
%

\Fref{fig:Lin} illustrates the set-up with Lin's method, shown in projection onto compactified Poincar{\'e} coordinates that represent $\R^3$ inside the unit sphere (not shown) centered at the origin $\mathbf{0}$. The plane in \fref{fig:Lin} is the common Lin section $\Sigma$ defined by $x = x_{\rm B} = 0$. Notice that the chosen orbit segment $\mathbf{u}$ intersects $\Sigma$ three times, that is, we choose to work with the third intersection of the trajectory from $\mathbf{0}$. Similarly, the chosen orbit segment $\mathbf{u}_{\rm B}$ intersects $\Sigma$ many times. The orbit segment $\mathbf{u}_{\rm B}$ in \fref{fig:Lin} was chosen so that its end point $\mathbf{u}_{\rm B}(1)$ in $\Sigma$ is sufficiently close to the end point $\mathbf{u}(0)$. The Lin space $Z \subset \Sigma$, which appears curved in the compactified Poincar{\'e} coordinates of \fref{fig:Lin}, remains fixed during the subsequent continuation of the BVP~\eqref{eq:BVPto0_vf}--\eqref{eq:BVPtoinfty_zero} and~\eqref{eq:BVPtoinftyL} in $\beta$. Panel~(b) shows the situation when the Lin gap $\eta$ has been closed and the connecting orbit found as the concatenation of $\mathbf{u}$ and $\mathbf{u}_{\rm B}$.  As is seen in \fref{fig:Lin}, the orbit segment $\mathbf{u}_{\rm B}$ intersects $\Sigma$ multiple times. We remark that, from a practical perspective, it is best to choose $\mathbf{u}_{\rm B}(1)$ close to $\mathbf{u}(0)$. On the other hand, choosing any of the earlier intersection points of $\mathbf{u}_{\rm B}$ in the numerical set-up will result in the same connecting orbit, provided that a zero of the Lin gap is found.

%
\begin{figure}[t!]
\begin{center}
  \includegraphics{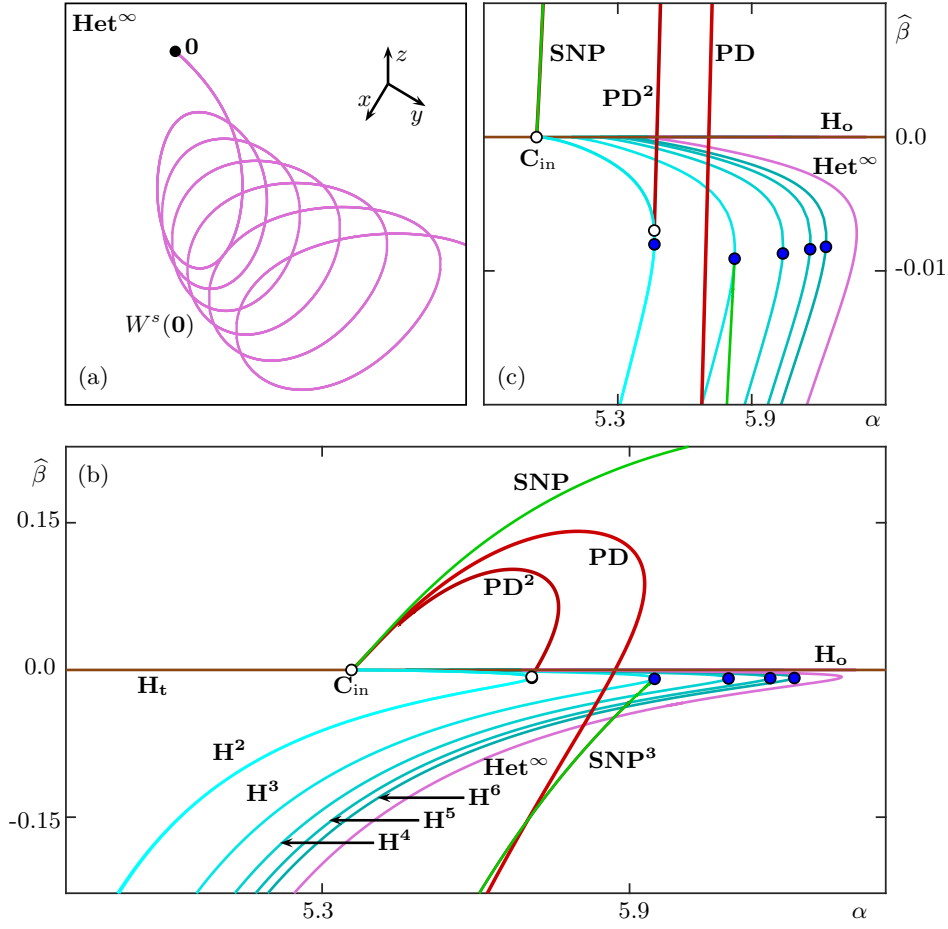}\label{fig:Formulation}
  \caption{Bifurcation diagram of system~\eqref{eq:WorkEquation} with the additional curve $\mathbf{Het^\infty}$ (magenta) of heteroclinic bifurcation involving the point $\mathbf{q}_\infty$ at infinity. Panel~(a) shows how $W^s(\mathbf{0})$ spirals towards infinity in the $(x,y,z)$-space to form the heteroclinic connection on $\mathbf{Het^\infty}$ for $\alpha = 5.3$ and $\beta=2.08874$; see \fref{fig:ConfigHom} for comparison. Panel~(b) shows the overall bifurcation diagram in the $(\alpha,\hat{\beta})$-plane and panel~(c) is an enlargement near the point $\mathbf{C}_{\rm in}$; see \fref{fig:bifDiagram} for details on the other bifurcation curves.} 
\end{center}
\end{figure}
%

As soon as a heteroclinic connection $\mathbf{Het^\infty}$  is detected as a zero of $\eta$, 
it can be continued with the BVP~\eqref{eq:BVPto0_vf}--\eqref{eq:BVPtoinfty_zero} and~\eqref{eq:BVPtoinftyL} in $\alpha$ and $\beta$, where $\theta_{\rm B} \in [0,\, 2 \pi]$, $T_0 > 0$, and $T_{\rm B} > 0$ are free parameters but $\eta = 0$ is now kept fixed. This continuation leads to the curve $\mathbf{Het^\infty}$ in the $(\alpha, \beta)$-plane that is shown in \fref{fig:Formulation} together with the other curves of the bifurcation diagram from \fref{fig:bifDiagram}. As panels~(a) and~(b) of \fref{fig:Formulation} show, the curve $\mathbf{Het^\infty}$ has the same general shape as the curves $\mathbf{H^n}$ of $n$-homoclinic bifurcation (shades of cyan) for $n = 2, 3, \dots, 6$: it also emanates from the codimension-two flip bifurcation point $\mathbf{C}_{\rm in}$, has monotonically decreasing $\hat{\beta}$ and has a fold for a very similar value of $\hat{\beta}$. Indeed, we conclude from \fref{fig:bifDiagram} that the curves $\mathbf{H^n}$ accumulate on the curve $\mathbf{Het^\infty}$ as $n$ tends to infinity. 

%
\begin{figure}[t!]
\begin{center}
  \includegraphics{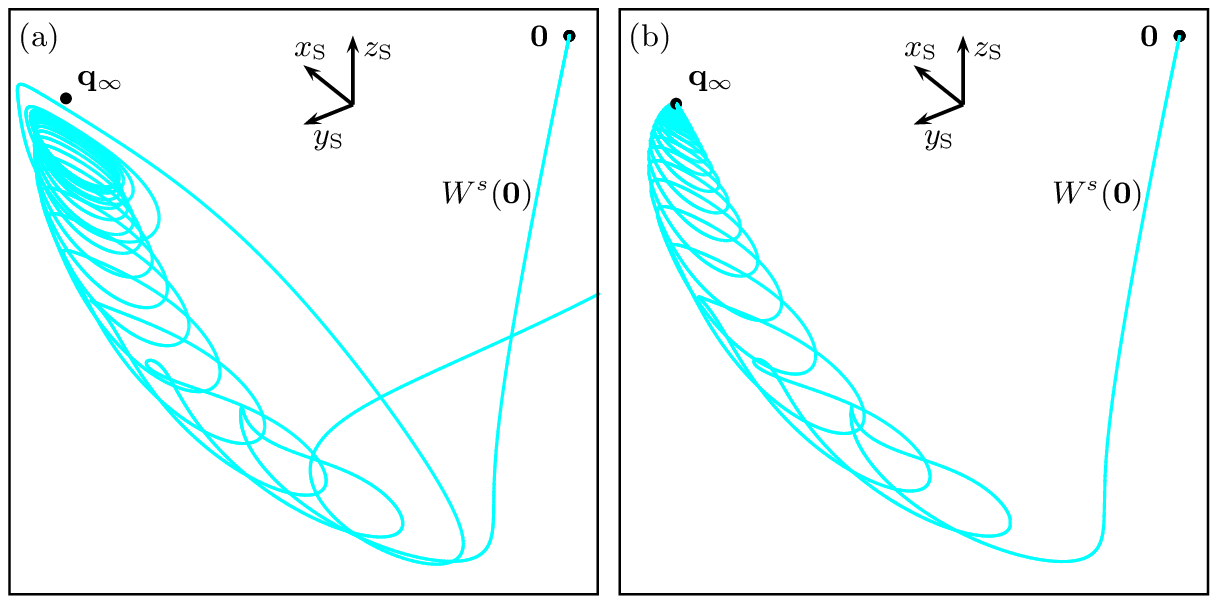}
  \caption{\label{fig:relativetoSc}
    To the left of the curve $\mathbf{Het^\infty}$ in the $(\alpha,\hat{\beta})$-plane, the stable manifold of $W^s(\mathbf{0})$ approaches, but does not connect to $\mathbf{q}_\infty$, because it lies outside $S_{\rm c}$ (a). To the right of $\mathbf{Het^\infty}$, it lies inside $S_{\rm c}$ and so connects to $\mathbf{q}_\infty$. The illustration in compactified Poincar{\'e} coordinates is for $\alpha = 5.3$ with $\beta = 2.9$ in panel~(a) and $\beta = 2.8$ in panel~(b).}
\end{center}
\end{figure}
%

Panel $\mathbf{Het^\infty}$ of \fref{fig:Formulation} illustrates that the heteroclinic connection from $\mathbf{0}$ to $\mathbf{q}_\infty$ is characterized by the one-dimensional manifold $W^s(\mathbf{0})$ spiraling away (in backward time) from $\mathbf{0}$ towards infinity to approach $\mathbf{q}_\infty$ along the cone/cylinder $S_{\rm c}$. Indeed, this is the limiting case between the two generic situations that are illustrated in \fref{fig:relativetoSc}. Either $W^s(\mathbf{0})$ lies outside $S_{\rm c}$ and does not reach $\mathbf{q}_\infty$, as in panel~(a), or it lies inside $S_{\rm c}$ and spirals onto $\mathbf{q}_\infty$, as in panel~(c). The former situation occurs to the left of the curve $\mathbf{Het^\infty}$ in the $(\alpha,\hat{\beta})$-plane of \fref{fig:Formulation}, while $W^s(\mathbf{0})$ connects generically to $\mathbf{q}_\infty$ to the right of $\mathbf{Het^\infty}$.

\section{BVP set-up for computing a generic connection from a saddle periodic orbit to infinity}
\label{sec:LinPeriodicInfty}
The Lin's method set-up from the previous section can be adapted to compute other types of connecting orbits to infinity. We demonstrate this here with the example of a heteroclinic connection from the orientable saddle periodic orbit $\Gamma_o$, which bifurcates from the curve $\mathbf{H_o}$ and exists for $\hat{\beta} > 0$, to the point $\mathbf{q}_\infty$. More specifically, we compute an orbit in the intersection set 
$W^s(\mathbf{\Gamma_o}) \cap S_{\rm c}$, which exists generically, because $W^s(\mathbf{\Gamma_o})$ and $S_{\rm c}$ are both two dimensional manifolds. As before, we concatenate two orbit segments: $\mathbf{u}$ from a common section $\Sigma$ to $\Gamma_o$ and $\mathbf{u}_{\rm B}$ from $\mathbf{q}_\infty$ to $\Sigma$, which are again found as solutions to the overall BVP~\eqref{eq:BVPto0_vf}--\eqref{eq:BVPtoinfty_zero} and~\eqref{eq:BVPtoinftyL}. The difference is that the vector $\mathbf{e}^s$ in boundary condition~\eqref{eq:BVPto0_zero} is now a vector in the stable Floquet bundle of $\Gamma_o$. The periodic orbit $\Gamma_o$ and its stable Floquet bundle can be computed and continued with the BVP set-up presented in~\cite{KraRie1}, yielding the vector $\mathbf{e}^s$ (for any value of the system parameters).

%
\begin{figure}[t!]
\begin{center}
  \includegraphics{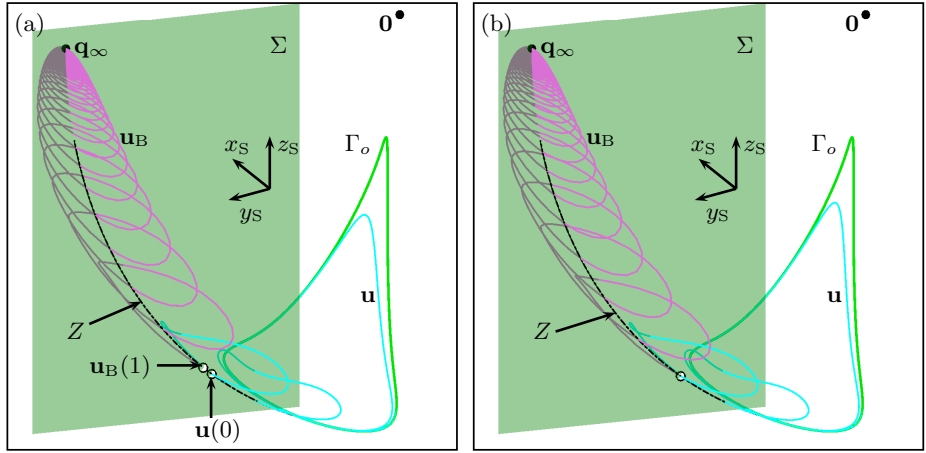}
  \caption{\label{fig:conPeriodic}
    Set-up with Lin's method to compute a connecting orbit from $\mathbf{q}_\infty$ to a saddle periodic orbit $\Gamma_o$ (green curve) with two orbit segments that meet in the common Lin section $\Sigma$ (green plane), illustrated in compactified Poincar{\'e} coordinates for $\alpha = 6.2$ and $\beta = 1.6$. Panel~(a) shows the initially chosen orbit segments $\mathbf{u}$ (cyan) to $\Gamma_o$ and $\mathbf{u}_{\rm B}$ (magenta) from $\mathbf{q}_\infty$ that define the Lin space $Z$ (which appears curved in this representation); note that the Lin gap $\eta$ is initially nonzero. Panel~(b) shows the situation where $\eta = 0$ and $\mathbf{u}$ and $\mathbf{u}_{\rm B}$ connect in $\Sigma$ to form the heteroclinic connection.}
  \end{center}
\end{figure}
%

A suitable initial orbit segment $\mathbf{u}$ is found by choosing and fixing $\delta_0$ and then, as before, continuing the initial value problem~\eqref{eq:BVPto0_vf} and~\eqref{eq:BVPto0_zero} in the integration time $T_0$ from $T_0 = 0$, while recording whenever condition~\eqref{eq:BVPto0_sigma} is satified. The initial orbit $\mathbf{u}_{\rm B}$ is found exactly as before, and the vector $\mathbf{\Psi}$, the Lin space $Z$ and the Lin gap $\eta$ are subsequently defined as in Section~\ref{sec:LinAtInfty}. The overall BVP~\eqref{eq:BVPto0_vf}--\eqref{eq:BVPtoinfty_zero} and~\eqref{eq:BVPtoinftyL} is then automatically satisfied and we use it to continue the two orbit segments $\mathbf{u}$ and $\mathbf{u}_{\rm B}$ to close the Lin gap $\eta$. Because the connecting orbit is generic, the continuation for this problem does not involve a system parameter, but uses the fact that the two-dimensional manifold $W^s(\mathbf{\Gamma_o})$ is a $\delta_0$-family of trajectories. Here, $\theta_{\rm B} \in [0,\, 2 \pi]$, $T_0 > 0$, $T_{\rm B} > 0$, $\eta \in \R$, and the parameter $\delta_0$ are free parameters.

\Fref{fig:conPeriodic} illustrates the set-up in compactified Poincar{\'e} coordinates; compare with \fref{fig:Lin}. Panel~(a) of \fref{fig:conPeriodic} shows the orientable periodic orbit $\Gamma_o$, the equilibrium $\mathbf{q}_{\infty}$, the section $\Sigma$ and the initially chosen orbit segments $\mathbf{u}$ and $\mathbf{u}_{\rm B}$ that define the Lin space $Z \subset \Sigma$. The Lin gap $\eta$ is then closed by continuation in $\delta_0$, yielding the connecting orbit as the concatenation of $\mathbf{u}$ and $\mathbf{u}_{\rm B}$ as shown in \fref{fig:conPeriodic}(b); note that the system parameters $\alpha$,  $\beta$ and $\gamma$ remain unchanged during this computation.

\section{Conclusions}
\label{sec:conclusions}

We studied a quadratic vector field, adapted from that of~\cite{Alg2019}, that exhibits a homoclinic flip bifurcation of the specific inward-twisted type $\mathbf{C}_{\rm in}$. We found that the two-parameter bifurcation diagram near this special point features an accumulation of curves of secondary $n$-homoclinic bifurcations. Numerical evidence that this phenomenon involves an increasing number of loops which move closer to infinity motivated us to set up a numerical scheme based on Lin's method to find the limiting behavior in the form of a heteroclinic connection to infinity. To this end, the orbit segment in the finite part of phase space was formulated in original coordinates, while the second orbit segment to infinity was defined in different coordinates near infinity. Both are then glued together along the Lin space in a section that is well-defined in both coordinate systems. Closing the Lin gap along the Lin space by continuation of the two coupled orbit segments yielded a first connecting orbit of codimension one between the origin and a point at infinity. A subsequent continuation gave the associated curve in the parameter plane, which was indeed found to act as the accumulation set for the curves of $n$-homoclinic orbits. 

Compared to previous uses of a Lin's method set-up to define suitable boundary value problems for finite connecting orbits, a novel element is the use of blown-up coordinate charts near infinity. Blow-up techniques for polynomial vector fields allow one to study equilibria and other invariant objects at and near infinity. When these are of saddle type in the geometric sense --- meaning that they have attracting and repelling directions, but need not be hyperbolic or even semi-hyperbolic --- the question arises how they interact with invariant objects in the finite part of phase space, such as equilibria and periodic orbits. Indeed, connections to infinity are a distinct possibility. As we showed, such heteroclinic phenomena involving infinity may provide important information regarding limits of finite global objects. 

Our Lin's method set-up is quite flexible and more widely applicable; this was demonstrated by computing a connecting orbit from a finite saddle periodic orbit to a point at infinity. Hence, it constitutes a new tool for the study of global properties of polynomial vector fields. The system studied here is a case in point, and its further bifurcation analysis is the subject of ongoing research. Note that this quadratic vector field is presently the only system that is known to exhibit a homoclinic flip bifurcation of the inward-twisted type of case $\mathbf{C}$; hence, it has the role of a model vector field for this specific bifurcation, much in the spirit of Sandstede's model~\cite{san4, Sandstede1997} which features effectively all other types of flip bifurcations. The investigation of the outward-twisted type in the latter model shows that a flip bifurcation of case $\mathbf{C}$ gives rise to a very complicated global bifurcation structure. In light of its different local structure, we expect to find a different, yet comparably complicated overall bifurcation structure in the wider vicinity of the orbit flip 
of the inward-twisted type of case $\mathbf{C}$. Moreover, homoclinic flip bifurcations of all cases have been identified as organizing centers in other vector fields from the literature, specifically in mathematical models of neurons~\cite{barrio2017, barrio2014, linaro2012}. Their global bifurcation structure may well involve heteroclinic bifurcations with infinity. Hence, we believe that the numerical approach for the identification and continuation of connecting orbits to infinity will have a role to play in their study.

\section*{Acknowledgments}
We thank Alejandro Rodr{\'\i}guez-Luis for alerting us to his paper~\cite{Alg2019} with Algaba, Dom{\'\i}nguez and Merino, and for suggesting that we investigate further the inward-twisted orbit flip bifurcation of case \textbf{C}.

%


\providecommand{\href}[2]{#2}
\providecommand{\arxiv}[1]{\href{http://arxiv.org/abs/#1}{arXiv:#1}}
\providecommand{\url}[1]{\texttt{#1}}
\providecommand{\urlprefix}{URL }

\end{document}